\documentclass[11pt]{article}

\usepackage[german,english]{babel}
\usepackage[cp850]{inputenc}
\usepackage{latexsym,graphicx}

\font\tenmath=msbm10 scaled 1200
\font\sevenmath=msbm7 scaled 1200
\font\fivemath=msbm5 scaled 1200

\newtheorem{Corollary}{Corollary}[section]

\newfam\mathfam \textfont\mathfam=\tenmath
\scriptfont\mathfam=\sevenmath \scriptscriptfont\mathfam=\fivemath
\def\math{\fam\mathfam}
\def\R{{\math R}}

\def\E \,{{\math E}}

\def\P{{\math P}}

\def \^#1{\if#1i{\accent"5E\i}\E \,lse{\accent"5E#1}\fi}
\def \E \,xo#1{{\bf Exercice #1} \mbox{}\\}

\def \E{{\math E}}

\def \ni{\noindent}

\newtheorem{Thm}{Theorem}

\author{
{\sc Harald Luschgy}\thanks{Universit\"at Trier, FB IV-Mathematik, D-54286
Trier, Germany.
E-mail: {\tt luschgy@uni-trier.de}} \quad {\sc  and}
\quad {\sc Gilles Pag\`es} \thanks{Laboratoire de Probabilit\'es et Mod\`eles
al\'eatoires, UMR~7599, Universit\'e Paris 6, case 188, 4,
pl. Jussieu, F-75252 Paris Cedex 5. E-mail:{\tt  gpa@ccr.jussieu.fr}}
}
\date{June 20, 06}
\title{\bf Moment estimates for L\'evy Processes}

\begin{document}

\maketitle
\begin{abstract}
For real L\'evy processes $(X_t)_{t \geq 0}$ having no Brownian component with
Blumenthal-Getoor index $\beta$,
the estimate $\E \, \sup_{s \leq t} | X_s - a_p s |^p \leq C_p t$ for
every $t \!\in [0,1]$ and suitable $a_p \!\in \R$ has been established by
Millar \cite{MILL} for $\beta < p \leq 2$ provided $X_1 \!\in L^p$. We
derive extensions of these estimates to the cases $p > 2$ and $p \leq
\beta$.
\end{abstract}

 
\bigskip
\ni {\em 2000 Mathematics Subject Classification:  60G51, 60G18.}


%
\noindent
\section{Introduction and results}
\setcounter{equation}{0}
\setcounter{Assumption}{0}
\setcounter{Theorem}{0}
\setcounter{Proposition}{0}
\setcounter{Corollary}{0}
\setcounter{Lemma}{0}
\setcounter{Definition}{0}
\setcounter{Remark}{0}
We investigate the $L^p$-norm (or quasi-norm) of the maximum process of real
L\'evy processes having no Brownian component. A (c\`adl\`ag) L\'evy process
$X=(X_t)_{t \geq 0}$ is characterized by its so-called local characteristics
in the L\'evy-Khintchine formula. They depend on the way the "big" jumps are
truncated. We will adopt in the following the convention that the truncation
occurs at size 1. So that
\begin{equation}
\E \, e^{i u X_t} = e^{-t \Psi (u)} \; \mbox{with} \; \Psi (u) = - i u a +
\frac{1}{2} \sigma2 u2 -\int (e^{i u x} - 1 - i u x \mbox{\bf 1 }_{ \{ | x |
\leq 1 \} } ) d \nu (x)
\end{equation}
where $u,\,a \!\in \R, \sigma2 \geq 0$ and $\nu$ is a measure on $\R$ such
that
$\nu (\{0\}) = 0$ and $\!\int x2 \wedge 1 d \nu (x) <+ \infty$.
The measure $\nu$ is called the L\'evy measure of $X$ and the quantities $(a,
\sigma2, \nu)$ are referred to as the characteristics of $X$.
One shows that for $p > 0, \E \, | X_1 |^p <+ \infty$ if and only if $\E \, |
X_t |^p <+ \infty$ for every $t \geq 0$ and this in turn is equivalent to
$\E \, \sup_{s \leq t} | X_s |^p <+ \infty$ for every $t \geq 0$. Furthermore,
\begin{equation}
\E \, | X_1 |^p <+ \infty \; \mbox{ if and only if } \;\int_{\{ x |  > 1 \} }
| x |^p d \nu (x) <+ \infty
\end{equation}
(see \cite{SATO}).
The index $\beta$ of the process $X$ introduced in \cite{BLUM} is defined by
\begin{equation}
\beta = \!\inf \{ p > 0 :\int_{ \{ | x |  \leq 1 \} } |  x |^p d \nu (x) <+
\infty \} .
\end{equation}
Necessarily, $\beta \!\in [0,2]$.

In the sequel we will assume that $\sigma2 = 0$, $i.e.$ that $X$ has no
Brownian component. Then the L\'evy-It\^o decomposition of $X$ reads
\begin{equation}
X_t = a t +\int^t_{0}\int_{\{ |  x | \leq 1 \} } x (\mu - \lambda \otimes
\nu)(ds, dx) +\int^t_0\int_{\{ | x |  > 1\} } x \mu (ds, dx)
\end{equation}
where $\lambda$ denotes the Lebesgue measure and $\mu$ is the Poisson random
measure on $\R_+ \times \R$
associated with the jumps of $X$ by
\[
\mu = \sum_{t \geq 0} \varepsilon_{(t, \triangle X_t)}  \mbox{\bf 1 }_{ \{
\triangle X_t \not= 0 \} } ,
\]
$\triangle X_t = X_t - X_{t-} , \triangle X_0 = 0$ (see \cite{JAWD} ,
\cite{SATO}).

\begin{Thm}
Let $(X_t)_{t \geq 0}$ be a L\'evy process with characteristics $(a, 0, \nu)$
and index $\beta$
such that
$\E \, | X_1 |^p <+ \infty$ for some
$p \!\in (\beta, \infty)$ or for $p = \beta$ provided
$\!\int_{\{ | x | \leq 1\} } | x |^\beta d \nu (x) <+ \infty$ and $\beta > 0$.
Then for every $t \geq 0$
\begin{eqnarray*}
\E \, \sup_{s \leq t} | Y_s |^p &\leq &C_p t \qquad \mbox{ if } \; p < 1 ,\\
\E \, \sup_{s \leq t} | X_s - s \,\E \, X_1 |^p &\leq& C_p t \qquad \mbox{ if
} \; 1 \leq p \leq 2
\end{eqnarray*}
where $Y_t = X_t - t (a -\int_{\{ | x |  \leq 1\} } x d \nu (x))$.
Furthermore, for every  $p > 2$
\[
\E \, \sup_{s \leq t}  | X_s |^p = O(t) \qquad \mbox{ as $t \rightarrow 0$}
\]
for a finite  real constant $C_p$.
\end{Thm}

If $X_1$ is symmetric one observes that $Y = X$ since the symmetry of $X_1$
implies $a = 0$ and the symmetry of $\nu$ (see  \cite{SATO}).
We emphasize that in view of the Kolmogorov criterion for continuous
modifications the above bounds are best possible as concerns powers of $t$. In
case $p > \beta $ and $p \leq 2$, these estimates are due to Millar
\cite{MILL}. However, the Laplace-transform approach in \cite {MILL} does not
work for $p > 2$. Our proof is based on the Burkholder-Davis-Gundy inequality.

For the case $p < \beta$ we need some assumptions on $X$.
Recall that a measurable function $\varphi \!:\! (0,c] \rightarrow (0,
\infty)\;(c > 0)$ is said to be regularly varying at zero with index $b \!\in
\R$ if, for
every
$t > 0$,
\[
\lim_{x \to 0} \frac{\varphi(tx)}{\varphi(x)} = t^b .
\]
This means that $\varphi(1/x)$ is regularly varying at infinity with index
$-b$. Slow variation corresponds to $b=0$.

\begin{Thm}
Let $(X_t)_{t \geq 0}$ be a L\'evy process with characteristics $(a, 0, \nu)$
and index $\beta$ such that
$\beta > 0$ and $\E \, | X_1 |^p <+ \infty$ for some $p \!\in (0, \beta)$.
Assume that the L\'evy measure satisfies
\begin{equation}
\exists\, c \!\in (0,1] , \mbox{\bf 1 }_{\{ 0 < | x | \leq c\}} \nu(dx) \leq
\varphi(| x | ) \mbox{\bf 1 }_{\{ 0 < | x |  \leq c \} } dx
\end{equation}
where $\varphi : (0,c] \rightarrow (0, \infty)$ is a regularly varying
function at zero of index $-(\beta+1)$. Let $l(x) = x^{\beta +1} \varphi(x)$
and assume that $l(1/x) , x \geq c$ is locally bounded. Let
$\underline{l}(x) = \underline{l}_\beta(x) = l (x^{1/\beta})$.

\smallskip
\noindent $(a)$ Assume $\beta > 1$. Then as $t \rightarrow 0$, for every $r
\!\in (\beta,2], q \!\in [p \vee 1, \beta)$,
\[
\E \, \sup_{s \leq t} | X_s |^p = O (t^{p/\beta} [ \underline{l} (t)^{p/r} +
\underline{l} (t)^{p/q}]) \quad \mbox{ if } \quad \beta < 2,
\]
\[
\E \, \sup_{s \leq t} | X_s |^p = O (t^{p/\beta} [ 1+ \underline{l} (t)^{p/q}
]) \quad \mbox{ if } \quad \beta = 2.
\]
If $\nu$ is symmetric then this holds for every $q \!\in[ p, \beta)$.

\smallskip
\noindent $(b)$   Assume $\beta < 1$. Then as $t \rightarrow 0$, for every $r
\!\in (\beta,1],\, q \!\in [p, \beta)$
\[
\E \, \sup_{s \leq t} | Y_s |^p = O ( t^{p/\beta} [\underline{l} (t)^{p/r} +
\underline{l} (t)^{p/q} ])
\]
where $Y_t = X_t - t (a -\int_{\{ | x | \leq 1\} } x d \nu (x))$ . If $\nu$ is
symmetric this holds for every $r \!\in (\beta,2]$.

\smallskip
\noindent  $(c)$ Assume $\beta = 1$ and $\nu$ is symmetric. Then as $t
\rightarrow 0$, for every $r \!\in (\beta, 2], q \!\in [p, \beta)$
\[
\E \, \sup_{s \leq t} | X_s - as |^p = O (t^{p/\beta} [ \underline{l}
(t)^{p/r} + \underline{l} (t)^{p/q}]) .
\]
\end{Thm}

It can be seen from strictly $\alpha$-stable L\'evy processes where $\beta =
\alpha$ that the above estimates are best possible as concerns powers of $t$.

Observe that condition (1.5) is satisfied for a broad class of L\'evy
processes. It implies that the tail function $t \mapsto \underline{\nu}(t) :=
\nu ([-t, t]^c), t > 0$ of the L\'evy measure is dominated, for $t \leq c$, by
$2\int^c_t \varphi(x) ds + \nu (| x |  > c)$, a regularly varying function at
zero with index $-\beta$, so that $\underline{\nu}(t) = O(t \varphi(t))$ as $t
\rightarrow 0$.

Important special cases are as follows.
\begin{Corollary} Assume the situation of Theorem 2 (with $\nu$ symmetric if
$\beta =1$)
and let $U$ denote any of the processes
$X,\,Y,\, (X_t - at)_{t \geq 0}$.

\smallskip
\noindent $(a)$ Assume that the slowly varying part $l$ of $\varphi$ is
decreasing and unbounded on
$(0, c] $ (e.g. $(-\log x)^a, a > 0)$. Then as $t \rightarrow 0$, for every
$\varepsilon \!\in (0, \beta)$,
\[
\E \, \sup_{ s \leq t} | U_s |^p = O (t^{p/\beta}
\underline{l}(t)^{p/(\beta-\varepsilon}).
\]

\smallskip
\noindent $(b)$ Assume that $l$ is increasing on $(0,c]$ satisfying $l(0+) = 0
$ (e.g. $(-\log x)^{-a}, a > 0, c < 1)$ and $\beta \!\in (0,2)$.
Then as $t \rightarrow 0$, for every $\varepsilon > 0$,
\[
\E \, \sup_{ s \leq t} | U_s |^p = O (t^{p/\beta} \underline{l}(t)^{p/(\beta +
\varepsilon}).
\]
\end{Corollary}

The remaining cases $p = \beta \!\in (0,2)$ if $\beta \not= 1$ and $p \leq 1$
if $\beta = 1$ are solved under the assumption
that the slowly varying part of the function $\varphi$ in (1.5) is constant.
\begin{Thm}
Let $(X_t)_{t \geq 0}$ be a L\'evy process with characteristics $(a, 0, \nu)$
and index $\beta$ such that $\beta \!\in (0,2)$ and
$\E \, | X_1 |^\beta <+ \infty$ if $\beta \not= 1$ and $\E \, | X_1 |^p <+
\infty$ for some $p \leq 1$ if $\beta = 1$.
Assume that the L\'evy measure satisfies
\begin{equation}
\exists\, c \!\in (0,1], \exists\, C \!\in (0, \infty) , \mbox{\bf 1 }_{\{ 0 <
| x | \leq c\}} \nu (dx) \leq \frac{C}{| x |^{\beta+1} }
\mbox{\bf 1 }_{\{0 < | x | \leq c\}} dx .
\end{equation}
Then as $t \rightarrow 0$
\begin{eqnarray*}
\E \, \sup_{s \leq t} | X_s |^\beta &=& O(t(-\log t)) \; \mbox{ if } \; \beta
>> 1, \\
\E \, \sup_{s \leq t} | Y_s |^\beta &=& O(t(-\log t)) \; \mbox{ if } \; \beta
< 1
\end{eqnarray*}
and
\[
\E \, \sup_{s \leq t} | X_s |^p = O((t(-\log t))^p) \quad \mbox{ if } \quad
\beta = 1, p \leq 1
\]
where the process $Y$ is defined as in Theorem 2.
\end{Thm}

The above estimates are optimal (see Section 3).

The paper is organized as follows. Section 2 is devoted to the proofs of
Theorems 1, 2 and 3. Section 3 contains a collection of examples.
\noindent
\section{Proofs}
\setcounter{equation}{0}
\setcounter{Assumption}{0}
\setcounter{Theorem}{0}
\setcounter{Proposition}{0}
\setcounter{Corollary}{0}
\setcounter{Lemma}{0}
\setcounter{Definition}{0}
\setcounter{Remark}{0}

We will extensively use the following compensation formula (see e.g.
\cite{JAWD})
\[
\E \,\int^t_0\int f(s, x) \mu (ds, dx) = \E \, \sum\limits_{s \leq t} f (s,
\Delta X_s) \mbox{\bf 1 }_{\{ \Delta X_s \not= 0\}}
=\int^t_0\int f(s, x) d \nu (x) ds
\]
where $f : \R_{+} \times \R \rightarrow \R_{+}$ is a Borel function. \\ \\
{\bf Proof of Theorem 1.}
Since $\E \, | X_1 |^p <+ \infty$ and
$p > \beta$ (or $p=\beta$ provided
$\!\int_{\{ | x |  \leq 1\} } | x |^\beta d \nu (x) <+ \infty$ and $\beta >
0$), it follows from (1.2) that
\[\int | x |^p d \nu (x) <+ \infty .
\]
CASE 1: $0 < p < 1$. In this case we have $\beta < 1$ and hence
$\!\int_{\{ x |  \leq 1\}} | x |  d \nu (x) <+ \infty$.
Consequently, $X$ a.s. has finite variation on finite intervals.
By (1.4),
\[
Y_t = X_t - t \left(a -\int_{\{ | x |  \leq 1 \} } x d \nu (x)\right)
=\int^t_0\int x \mu(ds, dx) = \sum_{s \leq t} \triangle X_s
\]
so that, using the elementary inequality $(u +  v)^p \leq u^p + v^p$,
\[
\sup_{s \leq t} | Y_s |^p \leq \left( \sum_{s \leq t} | \triangle X_s |
\right)^p \leq \sum_{s \leq t} |
\triangle X_s |^p =\int^t_0\int | x|^p \mu (ds, dx) .
\]
Consequently,
\[
\E \, \sup_{s \leq t} | Y_s |^p \leq t\int | x |^p d \nu (x)\; \mbox{for
every} \; t \geq 0.
\]
CASE 2: $1 \leq p \leq 2$. Introduce the martingale
\[
M_t  :=  X_t - t \E \, X_1  =  X_t - t \left(a+\int_{\{ | x| >1\} } x d \nu
(x)\right)
 =\int^t_0\int x (\mu - \lambda \otimes \nu) (ds, dx).
\]
It follows from the Burkholder-Davis-Gundy inequality (see \cite{KALL}) that
\[
\E \, \sup_{s \leq t} | M_s |^p \leq C \E \, [M]^{p/2}_t
\]
for some finite constant $C$. Since $p/2 \leq 1$, the quadratic variation
$[M]$ of $M$ satisfies
\[
[M]^{p/2}_t = \left( \sum_{s \leq t} | \triangle X_s |^2\right)^{p/2} \leq
\sum_{s \leq t} | \triangle X_s |^p
\]
so that
\[
\E \, \sup_{s \leq t} | M_s |^p \leq C t\int | x |^p d \nu (x) \; \mbox{for
every} \; t \geq 0.
\]
CASE 3: $p > 2$. One considers again the martingale L\'evy process $M_t = X_t
- t \E \, X_1$.
For $k \geq 1$ such that $2^k \leq p$, introduce the martingales
\[
N^{(k)}_t  :=\int^t_0\int | x |^{2^{k}} (\mu - \lambda \otimes \nu) (ds, dx)
 =  \sum_{s \leq t} |  \triangle X_s |^{2^{k}} - t\int | x |^{2^{k}} d \nu(x)
.
\]
Set $m := \max \{ k \geq 1 : 2^k < p \}$. Again by the Burkholder-Davis-Gundy
inequality
\begin{eqnarray*}
\E \, \sup_{s \leq t} | M_s |^p & \leq & C\, \E \, [M]^{p/2}_t \\
& = & C\, \E \, \left( t\int x2 d \nu (x) + N^{(1)}_t\right)^{p/2} \\
& = & \leq C \left( t^{p/2}\left (\int x2 d \nu (x)\right)^{p/2} + \E \, |
N^{(1)}_t |^{p/2}\right) \\
& \leq & C \,( t + \E \, | N^{(1)}_t |^{p/2})
\end{eqnarray*}
$\mbox{for every} \; t \!\in [0,1] $
where $C$ is a finite constant that may vary from line to line. Applying
successively the Burkholder-Davis-Gundy inequality to the martingales
$N^{(k)}$ and exponents $p/2^k > 1, 1 \leq k \leq m$, finally yields
\[
\E \, \sup_{s \leq t} | M_s |^p \leq C (t + \E \, [N^{(m)}]^{p/2^{m+1}}_t) \;
\mbox{ for every} \;
t \!\in [0,1] .
\]
Using $p \leq 2^{m+1}$, one gets
\[
[N^{(m)}]^{p/2 m+1}_t = \left(\sum_{s \leq t} | \triangle X_s
|^{2^{m+1}}\right)^{p/2^{m+1}} \leq \sum_{s \leq t} | \triangle X_s |^p
\]
so that
\[
\E \, \sup_{s \leq t} | M_s |^p \leq C\left(t + t\int | x |^p d \nu (x)\right)
\quad \mbox{ for every} \; t \!\in [0,1] .
\]
This implies $\E \, \sup_{s \leq t} | X_s |^p = O(t)$ as $t \rightarrow 0$.
\hspace*{\fill}{$\Box$}

\bigskip
\noindent
{\bf Proof of Theorems 2 and 3.}
Let $p \leq \beta$ and fix $c \!\in (0, 1 ]$.
Let $\nu_1 = \mbox{\bf 1 }_{\{ | x |  \leq c \} } \cdot \nu$ and
$\nu_2 = \mbox{\bf 1 }_{\{ | x |  > c \} } \cdot \nu$. Construct
L\'evy processes $X^{(1)}$ and
$X^{(2)}$ such that $X \stackrel{d}{=} X^{(1)} + X^{(2)}$ and $X^{(2)}$ is a
compound Poisson process with
L\'evy measure $\nu_2$. Then $\beta = \beta (X) = \beta (X^{(1)}),
\beta(X^{(2)}) = 0$,
$\E \, | X^{(1)}| ^q <+ \infty$ for every $q > 0$ and
$\E \, | X^{(2)}_1 |^p <+ \infty$. It follows e.g. from Theorem 1 that for
every $t \geq 0$,
\begin{equation}
\E \, \sup_{s \leq t} | X^{(2)}_s |^p \leq C_p t \qquad \mbox{ if } \; p < 1,
\end{equation}
\[
\E \, \sup_{s \leq t} | X^{(2)}-s\, \E \, X^{(2)}_1 |^p \leq C_p t \qquad
\mbox{ if } \; 1 \leq p \leq 2
\]
where $\E \, X^{(2)}_1 =\int x d \nu_2(x) =\int_{\{ | x | > c \} } x d \nu
(x)$ .

As concerns $X^{(1)}$, consider the martingale
\[
Z^{(1)}_t := X^{(1)}_t - t E X^{(1)}_1 = X^{(1)}_t - t \left(a -\int x
\mbox{\bf 1 }_{\{ c < | x | \leq 1\} } d \nu(x)\right) =\int^t_0\int x
(\mu_1 - \lambda \otimes \nu_1)(ds, dx)
\]
where $\mu_1$ denotes the Poisson random measure associated with the jumps of
$X^{(1)}$.
The starting idea is to part the 'small' and the 'big' jumps of $X^{(1)}$ in a
non homogeneous way with respect to the function
$s \mapsto s^{1/\beta}$. Indeed one may decompose $Z^{(1)}$ as follows
\[
Z^{(1)} = M + N
\]
where
\[
M_t :=\int^t_0\int x \mbox{\bf 1 }_{\{ | x |  \leq s^{1 / \beta} \} }
(\mu_1 - \lambda \otimes \nu_1 )( ds, dx)
\]
and
\[
N_t  :=\int^t_0\int x \mbox{\bf 1 }_{\{ | x | > s^{1/\beta}\}} (\mu_1 -
\lambda \otimes \nu_1)(ds, dx)
\]
are martingales. Observe that for every $q > 0$ and $t \geq 0$,
\begin{eqnarray*}\int^t_0\int | x |^q \mbox{\bf 1 }_{\{ | x | > s^{1/\beta}\}}
d \nu_1 (x) ds
& = &\int | x |^q (| x |^\beta \wedge t) d \nu_1 (x) \\
& \leq &\int_{\{ | x |  \leq c \} } | x |^{\beta+q} d \nu(x) <+ \infty .
\end{eqnarray*}
Consequently,
\[
N_t =\int^t_0\int x \mbox{\bf 1 }_{\{ | x | > s^{1/\beta}\}} d \mu_1 (s,x) -
\psi(t)
\]
where $\psi (t) :=\int^t_0\int x \mbox{\bf 1 }_{\{ | x | > s^{1/\beta}\}} d
\nu_1(x) ds$. Furthermore, for every $r > \beta$ or
$r = 2$ and $t \geq 0$
\begin{equation}\int^t_0\int | x |^r \mbox{\bf 1 }_{\{ | x |  \leq
s^{1/\beta}\}} d \nu_1 (x) ds \leq t\int_{\{ | x | \leq c \}} | x |^r d \nu
(x) <+ \infty .
\end{equation}
In the sequel let $C$ denote  a finite constant that may vary from line to
line.

We first claim that for every $t \geq 0, r \!\in (\beta,2] \cap [1,2]$ and for
$r = 2$,
\begin{equation}
\E \, \sup_{s \leq t} | M_s | \leq C(\int^t_0\int | x |^r \mbox{\bf 1 }_{\{ |
x | \leq s^{1/\beta} \} } d \nu_1 (x) ds )^{p/r} .
\end{equation}
In fact, it follows from the Burkholder-Davis-Gundy inequality and from
$p/r \leq 1, r/2 \leq 1$ that
\begin{eqnarray*}
\E \, \sup_{s \leq t} | M_s |^p & \leq & \left( \E \, \sup_{s \leq t} | M_s
|^r \right)^{p/r} \\
& \leq & C \left(\E \, [M]^{r/2}_t \right )^{p/r} \\
& = & C\left (\E \left( \sum_{s \leq t} | \triangle X^{(1)}_s |^2 \mbox{\bf 1
}_{\{ | \triangle X^{(1)}_s | \leq
s^{1/ \beta} \} } \right)^{r/2} \right)^{p/r} \\
& \leq & C\left(\E \, \sum_{s \leq t} | \triangle X^{(1)}_s |^r \mbox{\bf 1
}_{\{ | \triangle X^{(1)}_s | \leq s^{1/\beta} \} } \right)^{p/r} \\
& = & C\left(\!\int^t_0\int | x |^r \mbox{\bf 1 }_{\{ | x | \leq s^{1/\beta}
\} } d \nu_1 (x) ds\right)^{p/r} .
\end{eqnarray*}
Exactly as for $M$, one gets for every $t \geq 0$ and every $q \!\in [p,2]
\cap [1,2]$ that
\begin{equation}
\E \, \sup_{s \leq t} | N_s |^p \leq C\left(\int^t_0\int | x |^q \mbox{\bf 1
}_{\{ | x | > s^{1/\beta}\}}
d \nu_1 (x) ds\right)^{p/q} .
\end{equation}
If $\nu$ is symmetric then (2.4) holds for every $q \!\in [p,2]$ (which of
course provides additional information in case $p < 1$ only). Indeed, $\psi =
0$ by the symmetry of $\nu$ so that
\[
N_t =\int^t_0\int x \mbox{\bf 1 }_{\{ | x | > s^{1/\beta} \}} d \mu_1 (s,x)
\]
and for $q \!\in [p,1]$
\begin{eqnarray}
 \E \, \sup_{s \leq t}  \left|\int^s_0\int x \mbox{\bf 1 }_{\{ | x |  >
u^{1/\beta} \} }\mu_1 (d u, dx) \right|^p
&  \leq &   \left( \E \, \sup_{s \leq t} \left|\int^t_0\int x \mbox{\bf 1
}_{\{ | x |  > u^{1/\beta} \}} \mu_1 (du, dx) \right|^q \right)^{p/q}  \\
& \leq  & \left( \E \, \sum_{s \leq t} \left| \triangle X^{(1)}_s \right|^q
\mbox{\bf 1 }_{\{ | \triangle X^{(1)}_s | > s^{1/\beta} \}} \right)^{p/q}
\nonumber \\
& = & \left(\int^t_0\int | x |^q \mbox{\bf 1 }_{\{ | x |  > s^{1/\beta} \}} d
\nu_1 (x) ds \right)^{p/q} . \nonumber
\end{eqnarray}

In the case $\beta < 1$ we consider the process
\begin{eqnarray*}
Y^{(1)}_t & := & Z^{(1)}_t + t\int x d \nu_1 (x) = X^{(1)}_t - t\,\left(a
-\int_{\{ | x | \leq 1\}} xd \nu (x)\right) \\
& = & M_t + N_t + t\int x d \nu_1 (x) \\
& = &\int^t_0\int x \mbox{\bf 1 }_{\{ |  x | \leq s^{1/\beta}\}} \mu_1 (ds,
dx) +\int^1_0\int x
       \mbox{\bf 1 }_{\{ | x | > s^{1/\beta} \}} \mu_1 (ds, dx).
\end{eqnarray*}
Exactly as in (2.5) one shows that for $t \geq 0$ and $r \!\in (\beta, 1]$
\begin{equation}
 \E \,  \sup_{s \leq t} \left|\int^s_0\int x \mbox{\bf 1 }_{\{ | x | \leq
u^{1/\beta} \}} \mu_1 (du, dx) \right|^p
 \quad \leq  \left(\int^t_0\int | x |^r
       \mbox{\bf 1 }_{\{ | x | \leq s^{1/\beta} \}} d \nu_1 (x)
ds\right)^{p/r}.
\end{equation}

Combining (2.1) and (2.3) - (2.6) we obtain the following estimates. Let
\[
Z_t = X_t - t\left(a -\int x \mbox{\bf 1 }_{\{ c < | x | \leq 1\}} d
\nu(x)\right).
\]
CASE 1: $\beta \geq 1$ and $p < 1$. Then for every $t \geq 0, r \!\in (\beta,
2] \cup \{ 2 \}, q \!\in [1,2]$,
\begin{eqnarray}
\E \, \sup_{s \leq t} | Z_s |^p & \leq & C \left(t + (\int^t_0\int | x |^r
\mbox{\bf 1 }_{\{ | x |  \leq s^{1/\beta}\}}
d \nu_1 (x) ds \right)^{p/r} \nonumber \\
& & +  \left(\int^t_0\int | x |^q \mbox{\bf 1 }_{\{ | x |  > s^{1/\beta}\}}
d \nu_1 (x) ds )^{p/q} \right).
\end{eqnarray}
If $\nu$ is symmetric (2.7) is even valid for every $q \!\in [p,2]$. \\ \\
CASE 2: $\beta \geq 1$ and $p \geq 1$. Then for every $t \geq 0, r \!\in
(\beta, 2] \cup \{ 2 \}, q \!\in [p,2]$,
\begin{eqnarray}
\E \, \sup_{s \leq t} | X_s - s \,\E \, X_1  |^p & \leq & C \left(t +
\left(\int^t_0\int | x |^r \mbox{\bf 1 }_{\{ | x |  \leq s^{1/\beta}\}} d
\nu_1 (x) ds
\right)^{p/r}\right.
\nonumber \\ & & + \left.\left (\int^t_0 | x |^q \mbox{\bf 1 }_{\{ | x |  >
s^{1/\beta}\}}
d \nu_1 (x) ds \right)^{p/q} \right).
\end{eqnarray}
CASE 3: $\beta < 1$. Then for every $t \geq 0, r \!\in (\beta,1], q \!\in
[p,1]$
\begin{eqnarray}
\E \, \sup_{s \leq t} | Y_s |^p & \leq & C \left(t + \left(\int^t_0\int | x
|^r \mbox{\bf 1 }_{\{ | x |  \leq s^{1/\beta}\}} d \nu_1 (x) ds
\right)^{p/r}\right.
\nonumber
\\ & & +\left.  \left(\int^t_0 | x |^q \mbox{\bf 1 }_{\{ | x |  >
s^{1/\beta}\}}
d \nu_1 (x) ds \right)^{p/q} \right).
\end{eqnarray}
If $\nu$ is symmetric then $Y = Z = (X_t - at)_{t \geq 0}$ and (2.9) is valid
for every $r \!\in (\beta, 2], q \!\in [p,2]$. \\

Now we deduce Theorem 2. Assume $p \!\in (0, \beta)$ and (1.5). The constant
$c$ in the above decomposition of $X$
is specified by the constant from (1.5). Then one just needs to investigate
the integrals appearing in the right hand side of the inequalities (2.7) -
(2.10).
One observes that Theorem 1.5.11 in [1] yields for
$r > \beta$,
\[\int | x |^r \varphi (| x | ) \mbox{\bf 1 }_{\{ | x |  \leq s^{1/\beta}\}}
dx \sim \frac{2}{r-\beta}
s^{\frac{r}{\beta}-1} l(s^{1/\beta} ) \; \mbox{as} \; s \rightarrow 0
\]
which in turn implies that for small $t$,
\begin{eqnarray}\int^t_0\int  | x |^r  \mbox{\bf 1 }_{\{ | x |  \leq
s^{1/\beta}\}} d \nu_1(x) ds
& \leq &\int^t_0\int | x |^r \varphi(| x| ) \mbox{\bf 1 }_{\{ | x |  \leq
s^{1/\beta}\}} dx ds  \\
& \sim & \frac{2 \beta}{(r-\beta)r} t^{r/\beta} l(t^{1/\beta}) \; \mbox{as} \;
t \rightarrow 0. \nonumber
\end{eqnarray}
Similarly, for $0 < q < \beta$,
\begin{eqnarray}\int^t_0\int  | x |^q  \mbox{\bf 1 }_{\{ | x |  >
s^{1/\beta}\}} d \nu_1(x) ds
& \leq &\int^t_0\int | x |^q \varphi(| x| ) \mbox{\bf 1 }_{\{ | x |  >
s^{1/\beta}\}} dx ds \\
& \sim & \frac{2 \beta}{(\beta-q)q} t^{q/\beta} l(t^{1/\beta}) \; \mbox{as} \;
t \rightarrow 0. \nonumber
\end{eqnarray}
Using (2.2) for the case $\beta = 2$ and $t + t^p = o (t^{p/\beta}
\underline{l} (t)^\alpha)$ as
$t \rightarrow 0, \alpha > 0$, for the case $\beta > 1$ one derives Theorem 2.

As for Theorem 3, one just needs a suitable choice of $q$ in (2.7) - (2.9).
Note that by (1.6) for every
$\beta \!\in (0,2)$ and $t \leq c^\beta$,
\[
\begin{array}{lll}\displaystyle \int^t_0\int  | x |^\beta  \mbox{\bf 1 }_{\{ |
x |  > s^{1/\beta}\}} d \nu_1(x) ds
&  = & \displaystyle \int^t_0\int | x |^\beta \mbox{\bf 1 }_{\{ c \geq | x |

>> s^{1/\beta}\}} d \nu(x) ds  \\

& \leq & \displaystyle C\int^t_0\int | x |^{-1} \mbox{\bf 1 }_{\{ c \geq | x |

>> s^{1/\beta}\}} dx ds \\

&  = & \displaystyle  C t (-\log t)
\end{array}
\]
so that $q = \beta$ is the right choice. (This choice of $q$ is optimal.)
Since by (2.10), for $r \!\in (\beta,2] ( \not= \emptyset)$,
\[\int^t_0\int | x |^r \mbox{\bf 1 }_{\{ | x | \leq s^{1/\beta}\}} d \nu_1 (x)
ds = O (t^{r/\beta})
\]
the assertions follow from (2.7) - (2.9). \hfill{$\Box$}

\section{Examples}

Let $K_\nu$ denote the modified Bessel function of the third kind and
index $\nu$given by
\[
K_\nu (z) = \frac{1}{2}\int^\infty_0 \!\!u^{\nu-1}\exp(- \frac{z}{2} (u +
\frac{1}{u})) d u,\qquad z > 0 .
\]
$\bullet$ {\em The $\Gamma$-process} is a subordinator (increasing L\'evy
process) whose distribution $\P_{X_t}$ at time $t > 0$ is a
$\Gamma(1,t)$-distribution
\[
\P_{X_t} (dx) = \frac{1}{\Gamma(t)} x^{t-1} e^{-x} \mbox{\bf 1 }_{(0, \infty)}
(x) dx .
\]
The characteristics are given by
\[
\nu(dx) = \frac{1}{x} e^{-x} \mbox{\bf 1 }_{(0, \infty)} (x) ds
\]
and $a =\int^1_0 x d \nu (x)=1 - e^{-1}$ so that $\beta = 0$ and $Y = X$. It
follows from Theorem 1 that
\[
\E \, \sup_{s \leq t} X^p_s = \E \, X^p_t = O (t)
\]
for every $p > 0$. This is  clearly the true rate since
\[
\E \, X^p_t = \frac{\Gamma (p+t)}{\Gamma(t+1)} t \sim \Gamma(p) t \; \mbox{as}
\; t \rightarrow 0.
\]
$\bullet$ The {\em $\alpha$-stable L\'evy Processes} indexed by
$\alpha \!\in (0,2)$ have L\'evy measure
\[
\nu(dx) = \left( \frac{C_1}{x^{\alpha+1}} \mbox{\bf 1 }_{(0, \infty)} (x) +
\frac{C_2}{| x |^{\alpha+1}} \mbox{\bf 1 }_{(-\infty,0)} (x) \right) dx
\]
with $C_i \geq 0, C_1 + C_2 > 0$ so that $\E \, | X_1 |^p <+ \infty$
for
$p \!\in (0, \alpha), \E \, | X_1 |^\alpha = \infty$ and $\beta = \alpha$. It
follows from Theorems 2 and 3 that for $p \!\in (0, \alpha)$,
\begin{eqnarray*}
\E \, \sup_{s \leq t} | X_s |^p & = & O(t^{\,p/\alpha}) \quad \mbox{ if } \quad
\alpha > 1 , \\
\E \, \sup_{s \leq t} | Y_s |^p & = & O(t^{\,p/\alpha}) \quad \mbox{ if } \quad
\alpha < 1 , \\
\E \, \sup_{s \leq t} | X_s |^p & = & O((t\,(- \log t))^p) \quad \mbox{ if }
\quad \alpha = 1 .
\end{eqnarray*}
Here Theorem 3 gives the true rate provided $X$ is not strictly stable. In
fact, if $\alpha = 1$ the scaling property in this case says that $X_t
\displaystyle \stackrel{d}{=} tX_1+Ct\log t$ for some real constant
$C\neq 0$ (see~\cite{SATO}, p.87) so that for $p<1$
\[
 \E\,|X_t|^p = t^p \E\,|X_1+C\log t|^p \sim C^p t^p |\log t|^p \qquad
\mbox{ as } \quad t\to 0.
\]
%
%
%

Now assume that $X$ is {\em strictly $\alpha$-stable}.
If $\alpha < 1$, then $a =\int_{| x | \leq 1} x d \nu (x)$ and thus $Y = X$
and if $\alpha = 1$, then $\nu$ is symmetric (see
{\cite{SATO}). Consequently, by Theorem 2, for every $\alpha \!\in (0,2), p
\!\in (0, \alpha)$,
\[
\E \, \sup_{s \leq t} | X_s |^p = O (t^{\,p/\alpha}).
\]
In this case Theorem 2 provides the true rate since the self-similarity
property of strictly stable L\'evy processes implies
\[
\E \, \sup_{s \leq t} | X_s |^p = t^{\,p/\alpha} \E \, \sup_{s \leq 1} | X_s |^p
.
\]
$\bullet$ {\em Tempered stable processes} are subordinators with L\'evy
measure
\[
\nu(d x) = \frac{2^\alpha \cdot \alpha}{\Gamma (1 - \alpha)} x^{-(\alpha+1)}
\exp( -\frac{1}{2} \gamma^{1/\alpha} x) \mbox{\bf 1 }_{(0, \infty)} (x) dx
\]
and first characteristic $a =\int^1_0 x d \nu(x), \alpha \!\in (0,1), \gamma
\!>\! 0 $ (see \cite{SCHOU}) so that
$\beta \!=\! \alpha,\, Y \!=\!X$ and $\E \, X^p_1 <+ \infty$ for every $p
> 0$. The distribution of $X_t$ is not generally known. It follows from
Theorems 1,2 and 3 that
\begin{eqnarray*}
 \E \, X^p_t &=& O(t) \quad \mbox{ if }
\quad p > \alpha , \\
\E \, X^p_t & = & O(t^{\,p/\alpha}) \qquad \quad\mbox{ if } \quad p < \alpha \\
\E \, X^\alpha_t & = & O(t(- \log t)) \quad \mbox{ if } \quad p = \alpha .
\end{eqnarray*}

For $\alpha=1/2$, the process reduces to the {\em inverse Gaussian
process} whose ditribution $\P_{X_t}$ at time $t>0$ is given by
\[
\P_{X_t}(dx) = \frac{t}{\sqrt{2\pi}}x^{-3/2}  \exp\left(-\frac
12(\frac{t}{\sqrt{x}}-\gamma\sqrt{x})^2\right)\mbox{\bf
1}_{(0,\infty)}(x)dx.
\]
In this case all rates are the true rates. In fact, for $p>0$, %
\begin{eqnarray*}
\E X_t^p&=&  \frac{t}{\sqrt{2\pi}} e^{t\gamma}\int_0^\infty
x^{p-3/2}\exp\left(-\frac 12(\frac{t^2}{x}+\gamma^2 x)\right)dx\\
&=&  \frac{t}{\sqrt{2\pi}}e^{t\gamma}\left(\frac
1\gamma\right)^{p-3/2}t^{p-1/2}
\int_0^\infty y^{p-3/2} \exp\left(-\frac{t\gamma}{2}(\frac{1}{y}+\gamma^2
y)\right)dy \\ 
&=& \frac{2}{\sqrt{2\pi}}\gamma^{p-3/2}t^{p+1/2} e^{t\gamma}
K_{p-1/2}(t\gamma)
\end{eqnarray*}
and, as $z\to 0$,
\begin{eqnarray*}
K_{p-1/2}(z)  & \sim & \frac{C_p}{z^{p-1/2}} \qquad \mbox{ if } \quad p
>\frac 12, \\
K_{p-1/2}(z)  & \sim & \frac{C_p}{z^{1/2-p}}  \quad \quad\mbox{ if }
\quad p <
\frac 12 \\
K_0(z)& \sim & |\log z| 
\end{eqnarray*}
where $C_p=2^{p-3/2}\Gamma(p-1/2)$ if $p>1/2$ and $C_p =
2^{-p-1/2}\Gamma(\frac 12-p)$ if $p<1/2$.

\medskip
\noindent $\bullet$ {\em The Normal Inverse Gaussian (NIG)} process
was introduced by Barndorff-Nielsen and has been used in financial modeling
(see \cite{SCHOU}). The NIG process
is a L\'evy process with characteristics $(a, 0, \nu)$ where
\[
\begin{array}{rcl}
\nu(d x) & = & \displaystyle{ \frac{\delta \alpha}{\pi} \; \frac{\exp
(\gamma x) K_1 (\alpha | x | ) }{ | x | } d x , } \\
a & = & \displaystyle{\frac{2 \delta \alpha}{\pi} \;\int^1_0 \sin h (\gamma x)
K_1 (\alpha x) dx , }
\end{array}
\]
$\alpha > 0, \;\gamma \!\in (- \alpha , \alpha), \delta > 0$.
Since $K_1 (| z | ) \sim  | z |^{-1 } $ as $z \rightarrow 0$, the L\'evy
density behaves like
$\delta \pi^{-1} | x |^{-2}$ as $ x \rightarrow 0$ so that (1.6) is
satisfied with $\beta=1$.
One also checks that $\E \, | X_1 |^p <+ \infty $ for every $p > 0$. It
follows from Theorems 1 and 3 that, as $t \rightarrow 0$
\begin{eqnarray*}
\E \, \sup_{s \leq t} | X_s |^p & = & O(t) \quad \mbox{ if } \quad p > 1 , \\
\E \, \sup_{s \leq t} | X_s |^p & = & O((t (- \log t))^p) \quad \mbox{ if
}
\quad p \leq 1 .
\end{eqnarray*}
If $\gamma = 0$, then $\nu$ is symmetric and by Theorem 2,
\[
\E \, \sup_{s \leq t} | X_s |^p = O(t^p) \quad \mbox{ if } \quad p < 1.
\]
The distribution $\P_{X_t}$ at time $t > 0$ is given by
\[
\P_{X_t} (dx) = \frac{t \delta \alpha}{\pi} \exp ( t\, \delta
\sqrt{\alpha^2 - \gamma^2} + \gamma x)
\frac{ K_1 ( \alpha \sqrt{t^2 \delta^2 + x^2} ) }{ \sqrt{t^2 \delta^2 +
x^2} } dx
\]
so that Theorem 3 gives the true rate for $p = \beta = 1$ in the symmetric
case. In fact, assuming $\gamma = 0$, we get as $t \rightarrow 0$
\begin{eqnarray*}
\E \, | X_t | & = & \displaystyle{ \frac{2 t \delta \alpha}{\pi} e^{t \delta
\alpha}\int^\infty_0
\frac{x K_1 ( \alpha \sqrt{t^2 \delta ^2 + x^2} ) }{ \sqrt{t^2 \delta^2 +
x^2} } dx } \\
& = & \displaystyle{ \frac{2 t \delta \alpha }{\pi} e^{t \delta
\alpha}\int^\infty_{t \delta} K_1 (\alpha y) dy } \\
& \sim & \displaystyle{ \frac{2 \delta }{\pi} t\int^1_{t \delta}
\frac{1}{y} dy } \\
& \sim & \displaystyle{ \frac{2 \delta}{\pi} t ( - \log (t )) }.
\end{eqnarray*}
$\bullet$ {\em Hyperbolic L\'evy motions} have been applied to option pricing
in finance (see \cite{EBER}). These processes are L\'evy processes whose
distribution $\P_{X_1}$ at time $t = 1$ is a symmetric (centered) hyperbolic
distribution
\[
\P_{X_1}  (dx) = C \,\exp ( -\delta \sqrt{1 + (x/ \gamma)^2}) dx ,
\quad \gamma ,
\delta > 0.
\]
Hyperbolic L\'evy processes have characteristics $(0,0,\nu)$ and satisfy $\E
\, | X_1 |^p <+ \infty$ for every $p > 0$. In particular, they are
martingales. There (rather involved)
symmetric L\'evy measure has a Lebesgue density that behaves like $C x^{-2} $
as $x \rightarrow 0$ so that (1.6) is satisfied with $\beta = 1$.
Consequently, by Theorems 1,2 and 3, as $t \rightarrow 0$
\begin{eqnarray*}
\E \, \sup_{s \leq t} | X_s |^p & = & O(t) \quad \mbox{ if } \quad p > 1 , \\
\E \, \sup_{s \leq t} | X_s |^p & = & O(t^p) \quad \mbox{ if } \quad p < 1 ,
\\
\E \, \sup_{s \leq t} | X_s | & = & O(t\,(- \log t)) \quad \mbox{ if }
\quad p = 1 .
\end{eqnarray*}
$\bullet$ {\em Meixner processes} are L\'evy processes without Brownian
component and with L\'evy measure given by
\[
\nu(dx) = \frac{\delta e^{\gamma x}}{x \,\mbox{sinh} (\pi x)} dx ,\; \delta >
0,\; \gamma \!\in (- \pi, \pi)
\]
(see \cite{SCHOU}). The density behaves like $\delta / \pi x2$ as $x
\rightarrow 0$ so that (1.6) is satisfied with $\beta = 1$. Using (1.2) one
observes that
$\E \, | X_1 |^p <+ \infty$ for every $p > 0$. It follows from Theorems 1 and
3 that
\begin{eqnarray*}
\E \, \sup_{s \leq t} | X_s |^p & = & O(t) \quad \mbox{ if } \quad p > 1 , \\
\E \, \sup_{s \leq t} | X_s |^p & = & O((t\,(- \log t))^p) \quad \mbox{ if }
\quad p \leq 1 .
\end{eqnarray*}
If $\gamma = 0$, then $\nu$ is symmetric and hence Theorem 2 yields
\[
\E \, \sup_{s \leq t} | X_s  |^p = O(t^p) \quad \mbox{ if } \quad  p < 1 .
\]
\end{document}